\documentclass{article}

\usepackage{amsmath}
\usepackage{amsfonts}
\usepackage{amssymb}
\usepackage{url}
\usepackage{color}
\usepackage{amssymb}
\usepackage{enumerate}
\usepackage{epsfig}

\newtheorem{theorem}{Theorem}
\newtheorem{lemma}[theorem]{Lemma}

\newtheorem{definition}[theorem]{Definition}
\newtheorem{corollary}[theorem]{Corollary}
\newtheorem{example}[theorem]{Example}
\newtheorem{remark}[theorem]{Remark}

\newcommand{\bc}[1]{\begin{corollary}\label{#1}}
\newcommand{\ec}{\end{corollary}}

\newcommand{\E}{{\mathbb E}}
\newcommand{\Q}{{\mathbb Q}}
\newcommand{\Z}{{\mathbb Z}}

\newcommand{\C}{{\mathbb C}}
\newcommand{\R}{\mathbb R}

\newcommand{\LP}{{({\text{LP}})}}
\newcommand{\IP}{{({\text{IP}})}}
\newcommand{\GLP}{{{\cal C}}}
\newcommand{\GIP}{{{\cal G}}}

\newcommand{\boproof}{{\bf Proof.} }
\newcommand{\eoproof}{\hspace*{\fill} $\square$ \vspace{5pt}}
\newcommand{\Orthant}{\mathbb O}

\DeclareMathOperator{\supp}{supp}

\def \F {{{\cal F}}}
\def \H {{{\cal H}}}

\usepackage{enumerate}

\begin{document}
\setlength{\parindent}{0pt} \setlength{\parskip}{2ex plus 0.4ex
minus 0.4ex}

\title{A polynomial oracle-time algorithm for convex integer
minimization}
\author{{\bf Raymond~Hemmecke}\\University of Magdeburg, Germany\\
\\
{\bf Shmuel Onn}\\Technion Haifa, Israel\\
\\
{\bf Robert Weismantel}\\University of Magdeburg, Germany}

\date{}

\maketitle

\begin{abstract}
In this paper we consider the solution of certain convex integer
minimization problems via greedy augmentation procedures. We 
show that a greedy augmentation procedure that employs only directions
from certain Graver bases needs only polynomially many augmentation
steps to solve the given problem. We extend these results to convex
$N$-fold integer minimization problems and to convex $2$-stage stochastic
integer minimization problems. Finally, we present some applications
of convex $N$-fold integer minimization problems for which our
approach provides polynomial time solution algorithms.
\end{abstract}

\section{Introduction}

For an integer matrix $A\in\Z^{d\times n}$, we define the circuits
$\GLP(A)$ and the Graver basis $\GIP(A)$ as follows. Herein, an
integer vector $v\in\Z^n$ is called \emph{primitive} if all its
components are coprime, that is, $\gcd(v_1,\ldots,v_n)=1$.

\begin{definition}
Let $A\in\Z^{d\times n}$ and let $\Orthant_j$, $j=1,\ldots,2^n$ denote
the $2^n$ orthants of $\R^n$. Then the cones
\[
C_j:=\ker(A)\cap\Orthant_j=\{z\in\Orthant_j:Az=0\}
\]
are pointed rational polyhedral cones. Let $R_j$ and $H_j$ denote the
(unique) minimal sets of primitive integer vectors generating $C_j$ over
$\R_+$ and $C_j\cap\Z^n$ over $\Z_+$, respectively. Then we define
\[
\GLP(A):=\bigcup_{j=1}^{2^n} R_j\setminus\{0\}\;\;\;\text{and}\;\;\;
\GIP(A):=\bigcup_{j=1}^{2^n} H_j\setminus\{0\}
\]
to be the set $\GLP(A)$ of circuits of $A$ and the Graver basis
$\GIP(A)$ of $A$.
\end{definition}

\begin{remark}\rm
It is not hard to show that $\GLP(A)$ corresponds indeed to
all primitive support-minimal vectors in $\ker(A)$ \cite{Graver:75}.
\end{remark}

Already in 1975, Graver showed that $\GLP(A)$ and $\GIP(A)$ provide
optimality certificates for a large class of continuous and integer
linear programs, namely for
\[
\LP_{A,u,b,f}:\qquad \min\{f(z):Az=b, 0\leq z\leq u, z\in\R_+^n\},
\]
and
\[
\IP_{A,u,b,f}:\qquad \min\{f(z):Az=b, 0\leq z\leq u, z\in\Z_+^n\},
\]
where the linear objective function $f(x)=c^\intercal x$, the upper
bounds vector $u$, and the right-hand side vector $b$ are allowed to
be changed \cite{Graver:75}. A solution $z^0$ to
$\LP_{A,u,b,f}$ is optimal if and only if there are no
$g\in\GLP(A)$ and $\alpha\in\R_+$ such that $z^0+\alpha g$ is a
feasible solution to $\LP_{A,u,b,f}$ that has a smaller objective
function value $f(z^0+\alpha g)<f(z^0)$. Analogously, an integer
solution $z^0$ to $\IP_{A,u,b,f}$ is optimal if and only if there are
no $g\in\GIP(A)$ and $\alpha\in\Z_+$ such that $z^0+\alpha g$ is
a feasible solution to $\IP_{A,u,b,f}$ that has a smaller objective
function value $f(z^0+\alpha g)<f(z^0)$.

Thus, the directions from $\GLP(A)$ and $\GIP(A)$ allow a simple
augmentation procedure that iteratively improves a given feasible
solution to optimality. While this augmentation process has to
terminate for bounded IPs, it may show some zig-zagging behaviour,
even to non-optimal solutions for LPs \cite{Hemmecke:PSP}:

\begin{example}\rm Consider the problem
\[
\min\{z_1+z_2-z_3: 2z_1+z_2\leq 2,z_1+2z_2\leq 2,z_3\leq
1,(z_1,z_2,z_3)\in\R_{\geq 0}^3\}
\]
with optimal solution $(0,0,1)$. Introducing slack variables
$z_4,z_5,z_6$ we obtain the problem $\min\{c^\intercal z:
Az=(2,2,1)^\intercal, z\in \R_{\geq 0}^6\}$ with
$c^\intercal=(1,1,-1,0,0,0)$ and
\[
A=\left(\begin{array}{cccccc}
  2 & 1 & 0 & 1 & 0 & 0\\
  1 & 2 & 0 & 0 & 1 & 0\\
  0 & 0 & 1 & 0 & 0 & 1\\
 \end{array}\right).
\]
The vectors $(1,0,0,-2,-1,0)$, $(0,1,0,-1,-2,0)$, $(1,-2,0,0,3,0)$,
$(2,-1,0,-3,0,0)$, $(0,0,1,0,0,-1)$ together with their negatives
are the circuits of $A$. The improving directions are given by all
circuits $v$ for which $c^\intercal v>0$.

\begin{center}
\includegraphics[height=5cm]{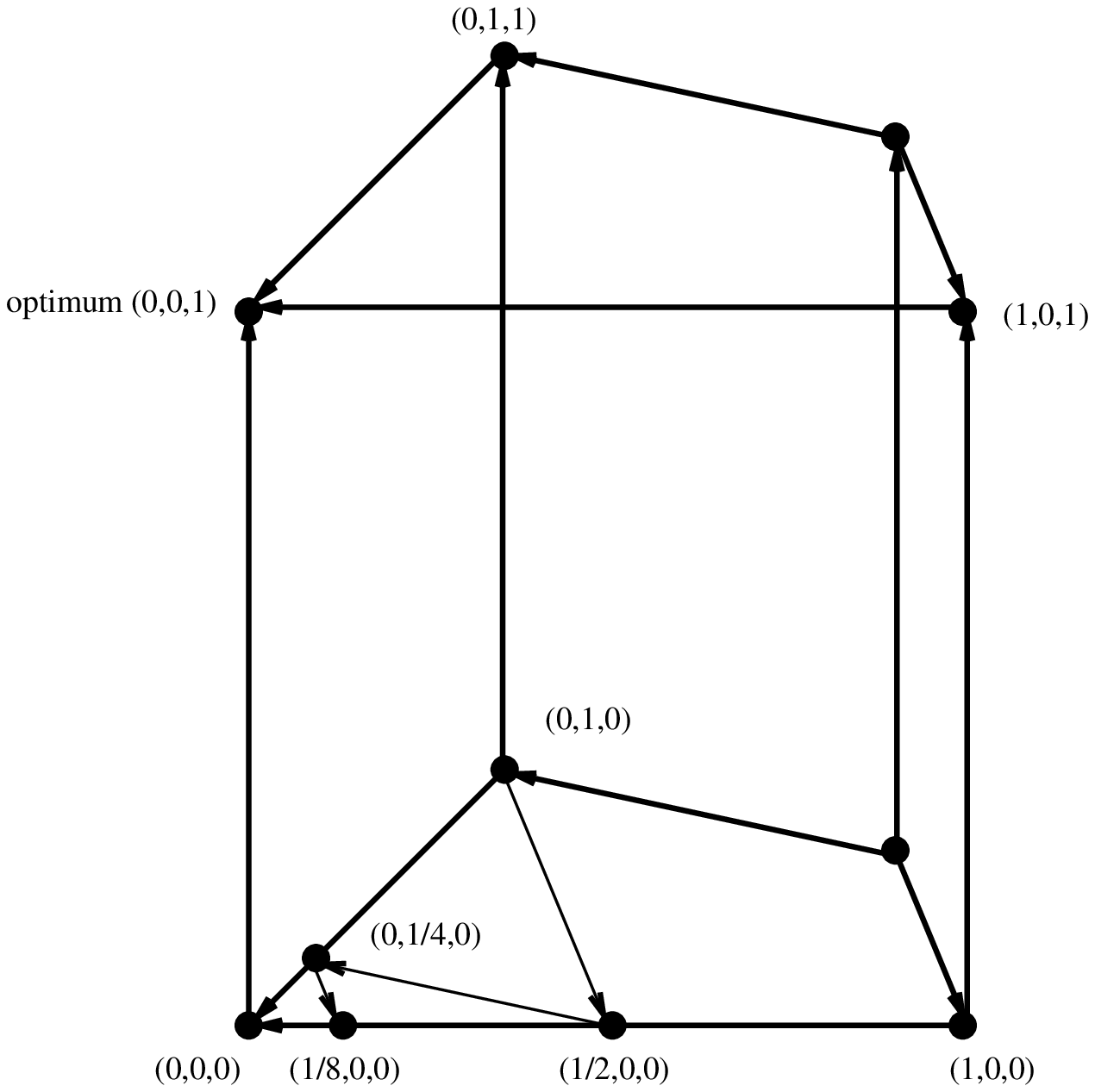}
\end{center}

Now start with the feasible solution $z_0=(0,1,0,1,0,1)$. Following
the directions $(0,1,0,-1,-2,0)$ and $(0,0,-1,0,0,1)$ as far
as possible, we immediately arrive at $(0,0,1,2,2,0)$ which
corresponds to the desired optimal solution $(0,0,1)$ of our
problem. However, alternatively choosing only the vectors
$(-1,2,0,0,-3,0)$ and $(2,-1,0,-3,0,0)$ as improving directions,
the augmentation process does not terminate. In our original space
$\R^3$, this corresponds to the sequence of movements
\[
(0,1,0) \rightarrow \left(\frac{1}{2},0,0\right)
      \rightarrow \left(0,\frac{1}{4},0\right)
      \rightarrow \left(\frac{1}{8},0,0\right)
      \rightarrow \left(0,\frac{1}{16},0\right)
      \rightarrow \ldots
\]
clearly shows the zig-zagging behaviour to the non-optimal point
$(0,0,0)$. \eoproof
\end{example}

Indeed, in order to avoid zig-zagging, certain conditions on the
selection of the potential augmenting circuits must be imposed.
As suggested in \cite{Hemmecke:PSP}, one can avoid such an undesired
convergence
\begin{itemize}
\item by first choosing an augmenting circuit direction freely, and
\item by then moving only along such circuit directions that do not
increase the objective value, that is $c^\intercal g\leq 0$, and which
introduce an additional zero component in the current feasible
solution, that is $\supp(z^0+\alpha g)\subsetneq\supp(z^0)$. After
$O(n)$ such steps, we have again reached a vertex and may perform a
free augmentation step if possible.
\end{itemize}

A natural question that arises is, whether there are strategies to
choose a direction from $\GLP(A)$ and $\GIP(A)$, respectively, to
augment any given feasible solution of $\LP_{A,u,b,f}$ or
$\IP_{A,u,b,f}$ to optimality in only polynomially many augmentation
steps. In this paper, we answer this question affirmatively. For this
let us introduce the notion of a \emph{greedy augmentation vector}.

\begin{definition}
Let ${\cal F}\subseteq\R^n$ be a set of feasible solutions,
$z_0\in{\cal F}$, $f:\R^n\rightarrow\R$ any objective function, and
let $S\subseteq\R^n$ be a (finite) set of directions. Then we call any
optimal solution to
\[
\min\{f(z_0+\alpha g):\alpha\in\R_+,g\in S,z_0+\alpha g\in{\cal F}\},
\]
a {\bf greedy augmentation vector} (from $S$ for $z_0$).
\end{definition}

\begin{theorem}\label{Theorem: Linear greedy augmentation algorithm}
Let $A\in\Z^{d\times n}$, $u\in\Q^n$, $b\in\Z^d$ and $c\in\Q^n$ be
given. Moreover, let $f(z)=c^\intercal z$. Then the following two
statements hold.
\begin{itemize}
  \item[(a)] Any feasible solution $z^0$ to $\LP_{A,u,b,f}$ can be
    augmented to an optimal solution of $\LP_{A,u,b,f}$ by iteratively
    applying the following greedy procedure:
    \begin{enumerate}
      \item Choose a greedy direction $\alpha g$ from $\GLP(A)$
    and set $z^0:=z^0+\alpha g$.
      \item[] If $\alpha g=0$, return $z_0$ as optimal solution.
      \item As long as it is possible, find a circuit direction
        $g\in\GLP(A)$ and $\alpha>0$ such that
    $z^0+\alpha g$ is feasible,
    $c^\intercal (z^0+\alpha g)\leq c^\intercal z^0$, and
    $\supp(z^0+\alpha g)\subsetneq\supp(z^0)$,
    and set $z^0:=z^0+\alpha g$.
      \item[] Go back to Step 1.
    \end{enumerate}
    The number of augmentation steps in this augmentation procedure is
    polynomially bounded in the encoding lengths of $A$, $u$, $b$,
    $c$, and $z^0$.
  \item[(b)] Any feasible solution $z^0$ to $\IP_{A,u,b,f}$ can be
    augmented to an optimal solution of $\IP_{A,u,b,f}$ by iteratively
    applying the following greedy procedure:
    \begin{enumerate}
      \item[]  Choose a greedy direction $\alpha g$ from $\GIP(A)$
    and set $z^0:=z^0+\alpha g$.
      \item[] If $\alpha g=0$, return $z_0$ as optimal solution.
    \end{enumerate}
    The number of augmentation steps in this augmentation procedure is
    polynomially bounded in the encoding lengths of $A$, $u$, $b$,
    $c$, and $z^0$.
  \end{itemize}
\end{theorem}

For our proof of Theorem \ref{Theorem: Linear greedy augmentation
algorithm} we refer to Section \ref{Subsection: Proof of linear greedy
algorithm}. Note that in \cite{DeLoera+Hemmecke+Onn+Weismantel} it was
shown that the Graver basis $\GIP(A)$ allows to design a polynomial time
augmentation procedure. This procedure makes use of the oracle
equivalence of so-called oriented augmentation and linear optimization
established in \cite{Schulz+Weismantel:99}. However, the choice of the 
Graver basis element that has to be used as a next augmenting vector
using the machanism of \cite{Schulz+Weismantel:99} is far more
technical than our simple greedy strategy suggested by Theorem
\ref{Theorem: Linear greedy augmentation algorithm}, Part (b).

In this paper, we generalize Part (b) of Theorem
\ref{Theorem: Linear greedy augmentation algorithm} to certain
$\Z$-convex objective functions. We say that a function
$g:\Z\rightarrow\Z$ is $\Z$-convex, if for all 
$x,y\in\Z$ and for all $0\leq\lambda\leq 1$ with $\lambda
x+(1-\lambda) y\in\Z$, the inequality $g(\lambda x+(1-\lambda) y)\geq
\lambda g(x)+(1-\lambda)g(y)$ holds. 
With this notion of $\Z$-convexity, we generalize Part (b) of Theorem
\ref{Theorem: Linear greedy augmentation algorithm} to
\emph{nonlinear convex} objectives of the form
$f(c^\intercal z,c_1^\intercal z,\ldots,c_s^\intercal z)$, where 
\begin{equation}\label{Equation: Separable function}
f(y_0,y_1,\ldots,y_s)=\sum_{i=1}^sf_i(y_i)+y_0
\end{equation}
is a separable $\Z$-convex function and where $c_0,\ldots,c_s\in\Z^n$
are given fixed vectors. In particular, each function
$f_i:\Z\rightarrow\Z$ is $\Z$-convex. When all $f_i\equiv 0$, we
recover linear integer optimization as a special case. To state our
result, let $C$ denote the $s\times n$ matrix with rows
$c_1,\ldots,c_s$ and let $\GIP(A,C)$ denote the Graver basis of 
$\left(\begin{smallmatrix}A & 0\\C & I_s\end{smallmatrix}\right)$
projected onto the first $n$ variables. As was shown in
\cite{Hemmecke:Z-convex,Murota+Saito+Weismantel}, this finite set
provides an improving direction for any non-optimal solution $z^0$ of
$\IP_{A,u,b,f}$.

\begin{theorem}\label{Theorem: Convex greedy augmentation algorithm}
Let $A\in\Z^{d\times n}$, $u\in\Z^n$, $b\in\Z^d$, $c\in\Q^n$,
$c_1,\ldots,c_s\in\Q^n$. Moreover, let $\bar{f}(z):=f(c^\intercal
z,c_1^\intercal z,\ldots,c_s^\intercal z)$, where $f$ denotes a
separable $\Z$-convex function as in (\ref{Equation: Separable
function}) given by a polynomial time comparison oracle which,
when queried on $x,y\in\Z^{s+1}$, decides whether $f(x)<f(y)$,
$f(x)=f(y)$, or $f(x)>f(y)$ holds in time polynomial in the encoding
lengths of $x$ and $y$. Moreover, let $H$ be an upper bound for the
difference of maximum and minimum value of $\bar{f}$ over the feasible
set $\{z:Az=b, 0\leq z\leq u, z\in\Z_+^n\}$ and assume that the
encoding length of $H$ is of polynomial size in the encoding lengths
of $A,u,b,c,c_1,\ldots,c_s$. Then the following statement holds. 

Any feasible solution $z^0$ to $\IP_{A,u,b,\bar{f}}$ can be augmented
to an optimal solution of $\IP_{A,u,b,\bar{f}}$ by iteratively
applying the following greedy procedure:
\begin{enumerate}
  \item[]  Choose a greedy direction $\alpha g$ from $\GIP(A,C)$
    and set $z^0:=z^0+\alpha g$.
  \item[] If $\alpha g=0$, return $z_0$ as optimal solution.
\end{enumerate}
The number of augmentation steps in this augmentation procedure is
polynomially bounded in the encoding lengths of $A$, $u$, $b$,
$c,c_1,\ldots,c_s$, and $z^0$.
\end{theorem}

For our proof of Theorem \ref{Theorem: Convex greedy augmentation
algorithm} we refer to Section \ref{Subsection: Proof of convex greedy
algorithm}. As a consequence to Theorem \ref{Theorem: Convex greedy
augmentation algorithm}, we construct in Sections \ref{Section: Convex
N-fold integer minimization} and \ref{Section: Convex 2-stage stochastic
integer minimization} polynomial time algorithms to solve convex
$N$-fold integer minimization problems and convex $2$-stage stochastic
integer minimization problems. In the first case, the Graver basis
under consideration is of polynomial size in the input data and hence
the greedy augmentation vector $\alpha g$ can be found in polynomial
time. In the second case, the Graver basis is usually of exponential
size in the input data. Despite this fact, the desired greedy
augmentation vector $\alpha g$ can be constructed in polynomial time,
if the $f_i$ are convex \emph{polynomial} functions. Finally, we
present some applications of convex $N$-fold integer minimization
problems for which our approach provides a polynomial time solution
algorithm. We conclude the paper with our proofs of Theorems
\ref{Theorem: Linear greedy augmentation algorithm} and \ref{Theorem:
Convex greedy augmentation algorithm}.

\section{$N$-fold convex integer minimization}
\label{Section: Convex N-fold integer minimization}

Let $A\in\Z^{d_a\times n}$, $B\in\Z^{d_b\times n}$, and
$c_1,\ldots,c_s\in\Z^n$ be fixed and consider the problem
\[
\min\left\{\sum_{i=1}^N f^{(i)}\left(x^{(i)}\right):\sum_{i=1}^N
Bx^{(i)}=b^{(0)}, Ax^{(i)}=b^{(i)}, 0\leq x^{(i)}\leq u^{(i)},
x^{(i)}\in\Z^n, i=1,\ldots,N\right\},
\]
where we have
\[
f^{(i)}(z):=\sum_{j=1}^sf^{(i)}_j\left(c_j^\intercal z\right)+{c^{(i)}}^\intercal z
\]
for given convex functions $f^{(i)}_j$ and vectors
$c^{(i)}\in\Z^n$, $i=1,\ldots,N$, $j=1,\ldots,s$. If we dropped the
coupling constraint $\sum_{i=1}^N Bx^{(i)}=b^{(0)}$, this optimization
problem would decompose into $N$ simpler convex problems
\[
\min\left\{f^{(i)}\left(x^{(i)}\right):Ax^{(i)}=b^{(i)}, 0\leq
x^{(i)}\leq u^{(i)}, x^{(i)}\in\Z^n\right\}, i=1,\ldots,N,
\]
which could be solved independently. Hence the name ``$N$-fold convex
integer program''.

\begin{definition}
The {\bf $N$-fold matrix of the ordered pair $A,B$} is the
following $(d_b+Nd_a)\times Nn$ matrix,
\[
[A,B]^{(N)}:= \left(
\begin{array}{ccccc}
  B & B & B & \cdots & B      \\
  A & 0 & 0 & \cdots & 0      \\
  0 & A & 0 & \cdots & 0      \\
  \vdots & \vdots & \ddots & \vdots & \vdots \\
  0  & 0 & 0      & \cdots & A      \\
\end{array}
\right).
\]
For any vector $x=(x^1,\ldots,x^N)$ with $x^i\in\Z^n$ for
$i=1,\ldots,N$, we call the number $|\{i:x^i\neq 0\}|$ of
nonzero building blocks $x^i\in\Z^n$ of $x$ the {\bf type} of $x$.
\end{definition}

In \cite{Hosten+Sullivant}, it was shown that there exists a constant
$g(A,B)$ such that for all $N$ the types of the Graver basis elements
in $\GIP([A,B]^{(N)})$ are bounded by $g(A,B)$. In
\cite{DeLoera+Hemmecke+Onn+Weismantel}, this was exploited to solve
\emph{linear} $N$-fold IP in polynomial time.

\begin{lemma}[Results from \cite{DeLoera+Hemmecke+Onn+Weismantel}]\mbox{}
\begin{itemize}
  \item For fixed matrices $A$ and $B$ the sizes of the Graver bases
  $\GIP([A,B]^{(N)})$ increase only polynomially in $N$.
  \item For any choice of the right-hand side vector $b$, an initial
    feasible solution $z_0$ can be constructed in time polynomial in
    $N$ and in the encoding length of $b$.
  \item For any linear objective function $c^\intercal z$, this
    solution $z_0$ can be augmented to optimality in time polynomial
    in $N$ and in the encoding lengths of $b$, $c$, $u$, and $z_0$.
\end{itemize}
\end{lemma}

Using Theorem \ref{Theorem: Convex greedy augmentation algorithm}, we
can now generalize this polynomial time algorithm to convex objectives
of the form above. Let us prepare the main result of this section by
showing that the encoding lengths of Graver bases from
\cite{Hemmecke:Z-convex,Murota+Saito+Weismantel} increase only
polynomially in $N$. For this, let $C$ denotes the $s\times n$ matrix
with rows $c_1,\ldots,c_s$.

\begin{lemma}\label{Lemma: ABC-Graver basis is of polynomial size}
Let the matrices $A\in\Z^{d_a\times n}$, $B\in\Z^{d_b\times n}$, and
$C\in\Z^{s\times n}$ be fixed. Then the encoding lengths of the Graver
bases of
\[
([A,B],C)^{(N)}:= \left(
\begin{array}{cccc|cccc}
  B & B & \cdots & B &&&& \\
  A &   &        &   &&&& \\
    & A &        &   &&&& \\
    &   & \ddots &   &&&& \\
    &   &        & A &&&&     \\
\hline
  C &   &        &   & I_s &     &        & \\
    & C &        &   &     & I_s &        & \\
    &   & \ddots &   &     &     & \vdots & \\
    &   &        & C &     &     &        & I_s \\
\end{array}
\right).
\]
increase only polynomially in $N$.
\end{lemma}

\boproof
The claim follows from the results in
\cite{DeLoera+Hemmecke+Onn+Weismantel} by rearranging the rows and
columns as follows
\[
([A,B],C)^{(N)}:= \left(
\begin{array}{cccccccc}
  B & 0    & B & 0      & \cdots & B & 0  \\
  A & 0    &   &        &        &   &    \\
  C & I_s  &   &        &        &   &    \\
    &      & A & 0      &        &   &    \\
    &      & C & I_s    &        &   &    \\
    &      &   & \ddots & \ddots &   &    \\
    &      &   &        &        & A & 0  \\
    &      &   &        &        & C & I_s \\
\end{array}
\right).
\]
This is the matrix of an $N$-fold IP with
$\bar{A}=\left(\begin{smallmatrix}A&0\\C&I_s\\\end{smallmatrix}\right)$
and with
$\bar{B}=\left(\begin{smallmatrix}B & 0\\\end{smallmatrix}\right)$.
Hence, the sizes and the encoding lengths of the Graver bases increase
only polynomially in $N$. \eoproof

Now that we have shown that the Graver basis is of polynomial size, we
can consider each Graver basis element $g$ independently and search
for the best $\alpha\in\Z_+$ such that $z_0+\alpha g$ is feasible and
has a smallest objective value. This can be done in polynomial time as
the following lemma shows.

\begin{lemma}\label{Lemma: best alpha can be found in polynomial time}
Let $f:\R\rightarrow\R$ be a convex function given by a comparision
oracle.
Then for any given numbers $l,u\in\Z$, the one-dimensional
minimization problem $\min\{f(\alpha):l\leq\alpha\leq u\}$ can be
solved by polynomially many calls to the comparision oracle.
\end{lemma}

\boproof
If the interval $[l,u]$ contains at most $2$ integers, return $l$ or
$u$ as the minimum, depending on the values of $f(l)$ and $f(u)$. If
the interval $[l,u]$ contains at least $3$ integers, consider the
integers $\lfloor (l+u)/2\rfloor-1$, $\lfloor (l+u)/2\rfloor$,
$\lfloor (l+u)/2\rfloor+1\in[l,u]$ and exploit convexity of $f$ to
bisect the interval $[l,u]$ as follows:

If $f(\lfloor (l+u)/2\rfloor-1)<f(\lfloor (l+u)/2\rfloor)$ holds, then
the minimum of $f$ must be attained in the interval $[l,\lfloor
(l+u)/2\rfloor]$. If, on the other hand, $f(\lfloor
(l+u)/2\rfloor)>f(\lfloor (l+u)/2\rfloor+1)$, then the minimum of $f$
must be attained in the interval $[\lfloor (l+u)/2\rfloor+1,u]$. If
none of the two holds, the minimum of $f$ is attained in the point
$\alpha=\lfloor (l+u)/2\rfloor$.

Clearly, after $O(\log (u-l))$ bisection steps, the minimization
problem is solved. \eoproof

The results in \cite{DeLoera+Hemmecke+Onn+Weismantel} together with
the previous two lemmas now immediately imply the main result of this
section.

\begin{theorem}\label{Theorem: Main theorem for convex N-fold IP}
Let $A$, $B$, $C$ be fixed integer matrices of appropriate dimensions.
Then the following holds. Moreover, let $f^{(i)}_j:\R\rightarrow\R$ be
convex functions mapping $\Z$ to $\Z$ given by polynomial time
evaluation oracles. Then the problem
\[
  \min\left\{\sum_{i=1}^N f^{(i)}\left(x^{(i)}\right):\sum_{i=1}^N
  Bx^{(i)}=b^{(0)}, Ax^{(i)}=b^{(i)}, 0\leq x^{(i)}\leq u^{(i)},
  x^{(i)}\in\Z^n, i=1,\ldots,N\right\},
\]
can be solved in time polynomial in the encoding length of the input
data.
\end{theorem}

\boproof Polynomial time construction of an initial feasible solution
from which we can start our augmentation process follows immediately
from the results in \cite{DeLoera+Hemmecke+Onn+Weismantel}.

To show that this feasible solution can be augmented to optimality in
polynomial time, we note that by Theorem \ref{Theorem: Convex greedy
augmentation algorithm} that only polynomially many greedy
augmentation steps are needed. By Lemma \ref{Lemma: ABC-Graver basis
is of polynomial size}, we only need to check polynomially many
directions $g$ to search for a greedy augmentation vector. But this
can be done in polynomial time by Lemma \ref{Lemma: best alpha can be
found in polynomial time}. \eoproof

\section{Convex $2$-stage stochastic integer minimization}
\label{Section: Convex 2-stage stochastic integer minimization}

Multistage stochastic integer programming has become an important
field of optimization, see
\cite{Birge+Louveaux:97,Louveaux+Schultz:03,Roemisch+Schultz:01} for
details. From a mathematical point of view, the data describing a
$2$-stage stochastic integer program is as follows. Let
$T\in\Z^{d\times m}$, $W\in\Z^{d\times n}$, 
$c_1,\ldots,c_s\in\Z^{m}$, $d_1,\ldots,d_s\in\Z^{n}$ be fixed, and
consider the problem
\[
\min\{\E_\omega(f^\omega(x,y)):Tx+Wy=b^\omega,0\leq x\leq u_x, 0\leq
y\leq u_y, x\in\Z^m, y\in\Z^n\},
\]
where $\omega$ is some probability distribution in a suitable
probability space and where $f$ is a convex function of the form
\[
f^\omega(x,y):=\sum_{j=1}^sf^\omega_j\left(c_j^\intercal
x+d_j^\intercal y\right)
\]
in which each $f^\omega_j:\R\rightarrow\R$ is a convex function.

Discretizing the probability distribution using $N$ scenarios, we
obtain the following convex integer minimization problem
{\small
\[
\min\left\{\sum_{i=1}^N f^{(i)}\left(x,y^{(i)}\right):
Tx+Wy^{(i)}=b^{(i)}, 0\leq x\leq u_x, 0\leq y^{(i)}\leq u^{(i)}_y,
x\in\Z^m, y^{(i)}\in\Z^n, i=1,\ldots,N\right\},
\]
}\noindent
where we have
\[
f^{(i)}(x,y):=\sum_{j=1}^sf^{(i)}_j\left(c_j^\intercal x+d_j^\intercal y\right)
\]
for given convex functions $f^{(i)}_j$. Note that fixing the
first-stage decision $x$ would decompose the optimization problem into
$N$ simpler convex problems
\[
\min\left\{f^{(i)}\left(x,y^{(i)}\right):Ax^{(i)}=b^{(i)},0\leq
y^{(i)}\leq u^{(i)}_y,y^{(i)}\in\Z^n\right\}, i=1,\ldots,N,
\]
which could be solved independently. However, the problem of finding
a first-stage decision $x$ with smallest overall costs would still
remain to be solved.

\begin{lemma}[Results from \cite{Hemmecke+Schultz:03}]\mbox{}
\begin{itemize}
  \item A vector $(v,w_1,\ldots,w_N)$ is in the kernel of the matrix
    \[
    \overline{[T,W]}^{(N)}:= \left(
    \begin{array}{ccccc}
      T & W & 0 & \cdots & 0      \\
      T & 0 & W & \cdots & 0      \\
      \vdots &   &   & \ddots &       \\
      T & 0 & 0 & \cdots & W      \\
    \end{array}
    \right)
    \]
    if and only if $(v,w_i)\in\ker\left(\overline{[T,W]}^{(1)}\right)$
    for all $i$, that is, if $Tv+Ww_i=0$ for all $i$.
  \item The Graver bases for the matrices $\overline{[T,W]}^{(N)}$
    decompose into a finite number of first-stage and second-stage
    building blocks that are independent on $N$.
  \item For any given \emph{linear} objective, any given right-hand side
    vector and any non-optimal feasible solution $z_0$, an improving
    vector to $z_0$ can be reconstructed from the building blocks in
    time linear in the number $N$ of scenarios.
\end{itemize}
\end{lemma}

Note that this finiteness result from \cite{Hemmecke+Schultz:03}
does not imply that the Graver basis of $\overline{[T,W]}^{(N)}$ is
of polynomial size in $N$. In fact, one can easily construct an
exponential size counter-example. Before we present the main result
of this section, let show that there exists a polynomial time
optimality certificate also for convex $2$-stage stochastic integer
minimization problems of the type above, if the matrices $T$ and $W$
are kept fix. For this, let $C$ denote the $s\times m$ matrix
with rows $c_1,\ldots,c_s$, and let $D$ denote the $s\times n$ matrix
with rows $d_1,\ldots,d_s$.

\begin{lemma}\label{Lemma: Polytime convex optimality certificate for SIP}
The Graver bases of the matrices
\[
\overline{[T,W,C,D]}^{(N)}:= \left(
\begin{array}{ccccccccc}
  T      & W &     &        &     &     &     &        &     \\
  T      &   &   W &        &     &     &     &        &     \\
  \vdots &   &     & \ddots &     &     &     &        &     \\
  T      &   &     &        &   W &     &     &        &     \\
  C      & D &     &        &     & I_s &     &        &     \\
  C      &   &   D &        &     &     & I_s &        &     \\
  \vdots &   &     & \ddots &     &     &     & \ddots &     \\
  C      &   &     &        &   D &     &     &        & I_s \\
\end{array}
\right)
\]
decompose into a finite number of first-stage and second-stage
building blocks that are independent on $N$.

For any given \emph{convex} objective, any given right-hand side
vector and any non-optimal feasible solution $z_0$, an improving
vector to $z_0$ can be reconstructed from the building blocks in
time linear in the number $N$ of scenarios.
\end{lemma}

\boproof To prove our first claim, we rearrange blocks within the
matrix $\overline{[T,W,C,D]}^{(N)}$ as follows:
\[
\left(
\begin{array}{cccccccc}
  T      & W &   0 &        &     &        &     &     \\
  C      & D & I_s &        &     &        &     &     \\
  T      &   &     &      W &   0 &        &     &     \\
  C      &   &     &      D & I_s &        &     &     \\
  \vdots &   &     &        &     & \ddots &     &     \\
  T      &   &     &        &     &        & W   &   0 \\
  C      &   &     &        &     &        & D   & I_s \\
\end{array}
\right)=\overline{\left[\left(\begin{array}{c}T\\C\end{array}\right),\left(\begin{array}{cc}W&0\\D&I_s\end{array}\right)\right]}^{(N)},
\]
which is the matrix of a $2$-stage stochastic integer program with
$N$ scenarios and fixed matrices
$\left(\begin{smallmatrix}T\\C\end{smallmatrix}\right)$ and
$\left(\begin{smallmatrix}W&0\\D&I_s\end{smallmatrix}\right)$.
Hence, its Graver basis consists out of a constant number of
building blocks independent on $N$. This proves the first claim.

To prove the second claim, note that the results from
\cite{Hemmecke:Z-convex,Murota+Saito+Weismantel} show that the
Graver basis of $\overline{[T,W,C,D]}^{(N)}$ projected down onto the
variables corresponding to $T$ and $W$ columns gives improving
directions for non-optimal solutions $z_0$ to
\[
\min\left\{\sum_{i=1}^N f^{(i)}\left(x,y^{(i)}\right):
Tx+Wy^{(i)}=b^{(i)}, 0\leq x\leq u_x, 0\leq y\leq u_y, x\in\Z^m,
y\in\Z^n, i=1,\ldots,N\right\}.
\]
Thus, these directions consist out of only a constant number of
building blocks independent on $N$. Let
$z=(x,y^{(1)},\ldots,y^{(N)})$ be a feasible solution and let
$g=(v,w^{(1)},\ldots,w^{(N)})$ be an augmenting vector formed out of
the constant number of first-stage and second-stage building blocks.
To be an improving direction, $g$ must satisfy the following
constraints:
\begin{itemize}
\item $T(x+v)+W(y^{(i)}+w^{(i)})=b^{(i)}$, $i=1,\ldots,N$,
\item $0\leq x+v\leq u_x$,
\item $0\leq y^{(i)}+w^{(i)}\leq u^{(i)}_y$, $i=1,\ldots,N$,
\item $\sum_{i=1}^N f^{(i)}\left(x+y,y^{(i)}+w^{(i)}\right)<\sum_{i=1}^N f^{(i)}\left(x,y^{(i)}\right)$.
\end{itemize}
For each of the finitely many first-stage building blocks perform
the following test: If $0\leq x+v\leq u_x$, try to find suitable
second-stage building blocks satisfying the remaining constraints,
which for fixed $v$ simplify to
\begin{itemize}
\item $Tv+Ww^{(i)}=0$, $i=1,\ldots,N$,
\item $0\leq y^{(i)}+w^{(i)}\leq u^{(i)}_y$, $i=1,\ldots,N$,
\item $\sum_{i=1}^N f^{(i)}\left(x+v,y^{(i)}+w^{(i)}\right)<\sum_{i=1}^N
f^{(i)}\left(x,y^{(i)}\right)$.
\end{itemize}
For fixed $v$, this problem decomposes into $N$ independent
minimization problems:
\[
\min\left\{f^{(i)}\left(x+v,y^{(i)}+w^{(i)}\right):Tv+Ww^{(i)}=0,0\leq
y^{(i)}+w^{(i)}\leq u^{(i)}_y\right\}, i=1,\ldots,N.
\]
If for those optimal values $\sum_{i=1}^N
f^{(i)}\left(x+v,y^{(i)}+w^{(i)}\right)<\sum_{i=1}^N
f^{(i)}\left(x,y^{(i)}\right)$ holds, we have found an improving
vector $g=(v,w^{(1)},\ldots,w^{(N)})$ for $z_0$. If one of these
minimization problems is infeasible or if $\sum_{i=1}^N
f^{(i)}\left(x+v,y^{(i)}+w^{(i)}\right)\geq\sum_{i=1}^N
f^{(i)}\left(x,y^{(i)}\right)$, then no augmenting vector for $z_0$
can be constructed using the first-stage building block $v$. If for
no first-stage building block $v$, an augmenting vector can be
constructed $z_0$ must be optimal. If there was an augmenting vector
for $z_0$ with some first-stage building block $v$, this vector or
even a better augmenting vector would have been constructed by the
procedure above when the first-stage building block $v$ was
considered. \eoproof

Note that the augmenting vector constructed in the proof of the
previous lemma need not be a Graver basis element (it may not be
minimal), but every Graver basis element could be constructed,
guaranteeing the optimality certificate. It remains to show how to
construct a \emph{greedy} augmentation vector from the building
blocks from the Graver basis. Note that the procedure in the
previous proof constructs an augmenting vector also for a
\emph{fixed} step length $\alpha$. To compute a greedy augmentation
vector, however, one has to allow $\alpha$ to vary. But then, the
minimization problem does not decompose into $N$ independent simpler
problems. It is this difficulty that enforces us to restrict the set
of possible convex functions.

\begin{definition}
We call a convex function $f:\R^{m+n}\rightarrow\R$ that maps
$\Z^{m+n}$ to $\Z$ {\bf splittable}, if for all fixed vectors $x\in\Z^m$,
$y,g_1,g_2\in\Z^n$, and for all finite intervals $[l,u]\subseteq\R$,
there exists polynomially many (in the encoding length of the problem
data) intervals $I_1,\ldots,I_r$ such that 
\begin{itemize}
  \item $[l,u]=\bigcup\limits_{i=1}^r I_r$,
  \item $I_i\cap I_j\cap\Z=\emptyset$ for all $1\leq i< j\leq r$, and
  \item for each $j=1,\ldots,r$, either 
    $f(x,y+\alpha g_1)\leq f(x,y+\alpha g_2)$ or
    $f(x,y+\alpha g_1)\geq f(x,y+\alpha g_2)$ holds for all 
    $\alpha\in I_j$.
\end{itemize}
\end{definition}

Note that convex polynomials of fixed maximal degree $k$ are splittable,
as $f(x,y+\alpha g_1)-f(x,y+\alpha g_2)$ switches its sign at most $k$
times. Hence each interval $[l,u]$ can be split into at most $k+1$
intervals with the desired property. With the notion of splittable
convex functions, we can now state and prove the main theorem of this
section. 

\begin{theorem}
Let $T$, $W$, $C$, $D$ be fixed integer matrices of appropriate
dimensions. Then the following holds.
\begin{itemize}
  \item[(a)] For any choice of the right-hand side vector $b$, an initial
    feasible solution $z_0$ to
{\small
    \[
      \min\left\{\sum_{i=1}^N f^{(i)}\left(x,y^{(i)}\right):
      Tx+Wy^{(i)}=b^{(i)}, 0\leq x\leq u_x, 0\leq y^{(i)}\leq u^{(i)}_y,
      x\in\Z^m, y^{(i)}\in\Z^n, i=1,\ldots,N\right\},
    \]
}\noindent
    can be constructed in time polynomial in $N$ and in the encoding
    length of the input data.
  \item[(b)] Then, for any choice of splittable convex functions
    $f^{(i)}$, this solution $z_0$ can be augmented to optimality in
    time polynomial in the encoding length of the input data. 
\end{itemize}
\end{theorem}

\boproof
Let us prove Part (b) first. This proof follows the main idea behind
the proof of Lemma \ref{Lemma: Polytime convex optimality certificate
for SIP}. Let $z=(x,y^{(1)},\ldots,y^{(N)})$ be a feasible solution
and let $g=(v,w^{(1)},\ldots,w^{(N)})$ be an augmenting vector formed
out of the constant number of first-stage and second-stage building
blocks. Again, for fixed $v$, we wish to consider each scenarios
independently. For this, note that the possible step length
$\alpha\in\Z_+$ is bounded from above by some polynomial size bound
$u_\alpha$, since our feasible region is a polytope. Since the convex
functions $f^{(i)}$ are splittable, we can for each scenario partition
the interval $[0,u_\alpha]$ into polynomially subintervals
$I_{i,1},\ldots,I_{i,r_i}$ such that for each interval $I_{i,j}$ there
is either no building block leading to a feasible solution or a
well-defined building block $w_{i,j}$ with $Tv+Ww_{i,j}=0$ and 
$0\leq y^{(i)}+\alpha w_{i,j}\leq u^{(i)}$ that minimizes
$f^{(i)}(x+v,y^{(i)}+\alpha w_{i,j})$ for all $\alpha\in I_{i,j}$. 
 
Taking the common refinement of all intervals $I_{i,j}$,
$i=1,\ldots,N$, $j=1,\ldots,r_i$, one obtains polynomially many
intervals $J_1,\ldots,J_t$, such that for each interval $J_i$ and for
all $\alpha\in J_i$, there is a well-defined building block for each
scenario minimizing the function value. For this fixed vector
$g=(v,w^{(1)},\ldots,w^{(N)})$ we then compute the best $\alpha\in
J_i$, and then compare these values $\sum_{i=1}^N
f^{(i)}\left(x+\alpha v,y^{(i)}+\alpha w^{(i)}\right)$ to find the
desired greedy augmentation vector. Applying Theorem \ref{Theorem:
Convex greedy augmentation algorithm}, this proves Part (b).

Finally, let us prove Part (a). For this, introduce nonnegative
integer slack-variables into the second-stages to obtain a linear IP
with problem matrix 
\[
\overline{[T,(W,I_d,-I_d)]}^{(N)}:= \left(
\begin{array}{ccccccccccc}
  T & W & I_d & -I_d & 0 &   0 &    0 & \cdots & 0 &   &      \\
  T & 0 &   0 &    0 & W & I_d & -I_d & \cdots & 0 &   &      \\
  \vdots & &  &      &   &     &      & \ddots &   &   &      \\
  T & 0 &   0 &    0 & 0 &   0 &    0 & \cdots & W & I_d & -I_d \\
\end{array}
\right)
\]
whose associated Graver basis is formed out of only constantly many
first- and second-stage building blocks. Using this extended
formulation, we may immediately write down a feasible solution. Using
only greedy directions from the Graver basis of
$\overline{[T,(W,I_d,-I_d)]}^{(N)}$, we can minimize the sum of all
slack-variables in polynomially many augmentation steps. Part (b) now
implies that an optimal solution to this extended problem can be found
in polynomial time. If all slack-variables are $0$, we have found a
feasible solution to our intial problem, otherwise the initial problem
is infeasible. \eoproof

Let us conclude with the remark that these polynomiality results for
convex $2$-stage stochastic integer minimization can be extended to
the multi-stage situation by applying the finiteness results from
\cite{Aschenbrenner+Hemmecke}.

\section{Some Applications}

Consider the following general nonlinear problems over an
arbitrary set $\F\subseteq\Z^n$ of feasible solutions:
\begin{description}

\item[Separable convex minimization:]
Find a feasible point $x\in\F$ minimizing a {\em separable} convex
cost function $f(x):=\sum_{i=1}^nf_i(x_i)$ with each $f_i$ a
univariate convex function. It generalizes standard linear
optimization with cost $f(x)=\sum_{i=1}^n c_i^\intercal x_i$ recovered
with $f_i(x_i):=c_i^\intercal x_i$ for some costs $c_i$. 

\item[Minimum $l_p$-distance:]
Find a feasible point $x\in\F$
minimizing the $l_p$-distance to a partially specified ``goal'' point
${\bar x}\in\Z^n$. More precisely, given $1\leq p\leq\infty$ and the
restriction ${\bar x}_I:=({\bar x}_i:i\in I)$ of ${\bar x}$ to a subset
$I\subseteq\{1,\dots,n\}$ of the coordinates, find $x\in\F$ minimizing
the $l_p$-distance 
$\|x_I-{\bar x}_I\|_p:=(\sum_{i\in I} |x_i-{\bar x}_i|^p)^{1\over p}$ for $1\leq p<\infty$
and $\|x_I-{\bar x}_I\|_\infty:=\max_{i\in I} |x_i-{\bar x}_i|$ for $p=\infty$.
\end{description}
Note that a common special case of the above is the
natural problem of $l_p$-norm minimization over $\F$,
$\min\{\|x\|_p:x\in\F\}$;
in particular, the $l_\infty$-norm minimization problem
is the min-max problem $\min\{\max_{i=1}^n |x_i|\,:\,x\in\F\}$.

In our discussion of $N$-fold systems below it will be convenient
to index the variable vector as $x=(x^1,\dots,x^N)$ with each block
indexed as $x^i=(x_{i,1},\dots,x_{i,n})$, $i=1,\dots,N$.

We have the following corollary of Theorem \ref{Theorem: Main theorem
for convex N-fold IP}, which will be used in the applications to follow.

\bc{separable_and_distance}
Let $A$ and $B$ be fixed integer matrices of compatible sizes.
Then there is an algorithm that, given any positive integer $N$,
right-hand sides $b^i$, and upper bound vectors $u^i$,
of suitable dimensions, solves the above problems over the
following set of integer points in an $N$-fold program
\begin{equation}\label{feasible}
\F\ =\ \{x=(x^1,\dots,x^N)\in\Z^{N\times n}\,:\,
\sum_{i=1}^N Bx^i=b^0,\ Ax^i=b^i,\ 0\leq x^i\leq u^i, \ i=1,\dots,N\}
\end{equation}
in time which is polynomial in $N$ and in the binary
encoding length of the rest of the input, as follows:
\begin{enumerate}
\item
For $i=1,\dots,N$ and $j=1,\dots,n$, let $f_{i,j}$ denote convex univariate
functions. Moreover, let $f(x):=\sum_{i=1}^N\sum_{j=1}^n f_{i,j}(x_{i,j})$
be given by a comparision oracle. Then the algorithm solves the
separable convex minimization problem 
$$\min\left\{\sum_{i=1}^N\sum_{j=1}^n f_{i,j}(x_{i,j})\,:\,x\in\F\right\}\ .$$

\item
Given any $I\subseteq\{1,\dots,N\}\times\{1,\dots,n\}$, any
partially specified integer point ${\bar x}_I:=(x_{i,j}:(i,j)\in I)$,
and any integer $1\leq p<\infty$ or $p=\infty$,
the algorithm solves the minimum $l_p$-distance problem
$$\min\left\{\|x_I-{\bar x}_I\|_p\,:\,x\in\F\right\}\ .$$
In particular, the algorithm solves the $l_p$-norm
minimization problem $\min\{\|x\|_p:x\in\F\}$.

\end{enumerate}
\ec

\boproof
Consider first the separable convex minimization problem. Then this is
just the special case of Theorem \ref{Theorem: Main theorem for convex
N-fold IP} with $c_j:={\bf 1}_j$ the standard $j$-th
unit vector in $\Z^n$ for $j=1,\dots,n$ and $c^i:=0$ in $\Z^n$ for $i=1,\dots,N$.
The objective function in Theorem \ref{Theorem: Main theorem for convex
N-fold IP} then becomes the desired objective,
$$\sum_{i=1}^N f^i(x^i)\ =\ \sum_{i=1}^N \sum_{j=1}^n
f_{i,j}(c_j^\intercal x^i)+c^ix^i\ =\
\sum_{i=1}^N\sum_{j=1}^n f_{i,j}(x_{i,j})\ .$$
Next consider the minimum $l_P$-distance problem. Consider first an integer
$1\leq p<\infty$. Then we can minimize the integer-valued $p$-th power
$\|x\|_p^p$ instead of the $l_p$-norm itself.
Define
$$f_{i,j}(x_{i,j})\ :=\
\left\{
  \begin{array}{ll}
    |x_{i,j}-{\bar x}_{i,j}|^p, & \hbox{if $(i,j)\in I$;} \\
    0, & \hbox{otherwsie.}
  \end{array}
\right.
$$
With these $f_{i,j}$, the objective in the separable convex minimization
becomes the desired objective,
$$\sum_{i=1}^N\sum_{j=1}^n f_{i,j}(x_{i,j})\ =\
\sum_{(i,j)\in I}|x_{i,j}-{\bar x}_{i,j}|^p\ =\ |x_I-{\bar x}_I|_p^p\ .$$
Next, consider the case $p=\infty$. Let $w:=\max\{|u_{i,j}|:i=1,\dots,N,j=1,\dots,n\}$
be the maximum upper bound on any variable.
We may assume $w>0$ else $\F\subseteq \{0\}$ and the integer program is trivial.
Choose a positive integer $q$ satisfying $q \log(1+(2w)^{-1})>\log (Nn)$.
Now solve the minimum $l_q$-distance problem and let $x^*\in\F$
be an optimal solution. We claim that $x^*$ also minimizes the
$l_\infty$-distance to $\bar x$. Consider any $x\in\F$. By standard
inequalities between the $l_\infty$ and $l_q$ norms,
$$\|x^*_I-{\bar x}_I\|_\infty\ \leq\ \|x^*_I-{\bar x}_I\|_q \ \leq\
\|x_I-{\bar x}_I\|_q\ \leq\ (Nn)^{1\over q} \|x_I-{\bar x}_I\|_\infty\ .$$
Therefore
$$\|x^*_I-{\bar x}_I\|_\infty-\|x_I-{\bar x}_I\|_\infty\ \leq\
((Nn)^{1\over q}-1)\|x_I-{\bar x}_I\|_\infty\ \leq\ ((Nn)^{1\over q}-1)2w\ <\ 1$$
where the last inequality holds by the choice of $q$.
Since $\|x^*_I-{\bar x}_I\|_\infty$ and $\|x_I-{\bar x}_I\|_\infty$
are integers we find that indeed
$\|x^*_I-{\bar x}_I\|_\infty\leq \|x_I-{\bar x}_I\|_\infty$
holds for all $x\in\F$ and the claim follows.
\eoproof

\subsection{Congestion-avoiding (multi-way) transportation and routing}

The classical (discrete) transportation problem is the following.
We wish to transport commodities (in containers or bins) on a traffic network
(by land, sea or air), or route information (in packets) on a communication network,
from $n$ suppliers to $N$ customers. The demand by customer $i$ is $d_i$ units
and the supply from supplier $j$ is $s_j$ units. We need to determine
the number $x_{i,j}$ of units to transport to customer $i$ from supplier $j$
on channel $i\leftarrow j$ subject to supply-demand requirements
and upper bounds $x_{i,j}\leq u_{i,j}$ on channel capacity so as to minimize
total delay or cost. The classical approach assumes a channel
cost $c_{i,j}$ per unit flow, resulting in linear total cost
$\sum_{i=1}^N\sum_{j=1}^n c_{i,j} x_{i,j}$.
But due to channel congestion when subject to heavy traffic or heavy
communication load, the transportation delay or cost on a channel are
actually a nonlinear convex function of the flow over it, such as
$f_{i,j}(x_{i,j})=c_{i,j}|x_{i,j}|^{\alpha_{i,j}}$ for suitable
$\alpha_{i,j}>1$, resulting in nonlinear total cost $\sum_{i,j}
f_{i,j}(x_{i,j})$, which is much harder to minimize.

It is often natural that the number of suppliers is small and fixed while
the number of customers is very large. Then the transportation problem is
an $N$-fold integer programming problem. To see this, index the variable
vector as $x=(x^1,\dots,x^N)$ with $x^i=(x_{i,1},\dots,x_{i,n})$
and likewise for the upper bound vector.
Let $b^i:=d_i$ for $i=1,\dots,N$ and let $b^0:=(s_1,\dots,s_n)$.
Finally, let $A=(1,\dots,1)$ be the $1\times n$ matrix with all entries equal
to $1$ and let $B$ be the $n \times n$ identity matrix. Then the $N$-fold
constraints $Ax^i=b^i$, $i=1,\dots,N$ and $B(\sum_{i=1}^N x^i)=b^0$
represent, respectively the demand and supply constraints. The feasible set
in (\ref{feasible}) then consists of the feasible transportations and
the solution of the congestion-avoiding transportation problem is
provided by Corollary \ref{separable_and_distance} part 1. So we have:

\bc{congestion-avoiding}
Fix the number of suppliers and let $f_{i,j}$, $i=1,\dots,N$,
$j=1,\dots,n$, denote convex univariate functions. Moreover, let 
$f(x):=\sum_{i=1}^N\sum_{j=1}^n f_{i,j}(x_{i,j})$ be given by a
comparision oracle. Then the congestion-avoiding transportation
problem can be solved in polynomial time.
\ec

This result can be extended to multi-way (high-dimensional) transportation problems
as well. In the {\em $3$-way line-sum transportation problem}, the set of feasible
solutions consists of all nonnegative integer $L\times M\times N$ arrays
with specified line-sums and upper bound (capacity) constraints,
\begin{equation}\label{feasible3}
\F\ :=\ \{x\in\Z^{L\times M\times N} \ :\ \sum_i x_{i,j,k}=r_{j,k}
\,,\ \sum_j x_{i,j,k}=s_{i,k} \,,\ \sum_k x_{i,j,k}=t_{i,j}
\,,\ 0\leq x_{i,j,k}\leq u_{i,j,k}\,\}\ .
\end{equation}
If at least two of the array-size parameters $L,M,N$ are variable then
even the classical linear optimization problem over $\F$ is NP-hard
\cite{DO1}. In fact, remarkably, {\em every} integer program {\em is}
a $3\times M\times N$ transportation program for some $M$ and $N$
\cite{DO2}. But when both $L$ and $M$ are relatively small and fixed,
the resulting problem over ``long'' arrays, with a large and variable
number $N$ of layers, is again an $N$-fold program. To see this, index
the variable array as $x=(x^1,\dots,x^N)$ with
$x^i=(x_{1,1,i},\dots,x_{L,M,i})$ and likewise for the upper bound
vector. Let $A$ be the $(L+M)\times LM$ incidence matrix of the
complete bipartite graph $K_{L,M}$ and let $B$ be the $LM \times LM$
identity matrix. Finally, suitably define the right-hand 
side vectors $b^h$, $h=0,\dots,N$ in terms of the given line sums $r_{j,k}$,
$s_{i,k}$, and $t_{i,j}$. Then the $n$-fold constraint $B(\sum_{h=1}^N
x^h)=b^0$ represents the line-sum constraints where summation over
layers occurs, whereas $Ax^h=b^h$, $h=1,\dots,N$, represent the
line-sum constraints where summations are within a single layer at a
time. Then we can minimize in polynomial time any separable convex
cost function $\sum_{i=1}^L\sum_{j=1}^M\sum_{k=1}^N
f_{i,j,k}(x_{i,j,k})$ over the set of feasible transportations $\F$ in
(\ref{feasible3}). So we have: 

\bc{congestion-avoiding-3}
Fix any $L$ and $M$ and let $f_{i,j,k}$, $i=1,\dots,L$,
$j=1,\dots,M$, $k=1,\dots,N$, denote convex univariate
functions. Moreover, let $f(x):=\sum_{i=1}^L\sum_{j=1}^M\sum_{k=1}^N
f_{i,j,k}(x_{i,j,k})$ be given by a comparision oracle. Then the
congestion-avoiding $3$-way transportation problem can be solved in
polynomial time.
\ec

Even more generally, this result holds for ``long'' $d$-way transportations
of any fixed dimension $d$ and for any {\em hierarchical} sum constraints,
see Section \ref{Subsection: Hierarchically-constrained multi-way
arrays} below for the precise definitions. 

\subsection{Error-correcting codes} 

Linear-algebraic error correcting codes generalize the ``check-sum''
idea as follows: a message to be communicated on a noisy channel is
arranged in a vector $x$. To allow for error correction, several sums
of subsets of entries of $x$ are communicated as well. Multi-way
tables provide an appealing way of organizing the check-sum
protocol. The sender arranges the message in a multi-way
$M_1\times\cdots\times M_d$ array $x$ and sends it along with the sums
of some of its lower dimensional sub-arrays (margins). The receiver
obtains an array $\bar x$ 
with some entries distorted on the way; it then finds an array $\hat x$ having
the specified check-sums (margins), that is $l_p$-closest to the received
distorted array $\bar x$, and declares it as the retrieved message.
For instance, when working over the $\{0,1\}$ alphabet, the useful
{\em Hamming distance} is precisely the $l_1$-distance. Note that the
check-sums might be distorted as well; to overcome this difficulty,
we determine ahead of time an upper bound $U$ on all possible check-sums,
and make it a fixed part of the communication protocol; then we blow each array
to size $(M_1+1)\times\cdots\times (M_d+1)$, and fill in the new ``slack''
entries so as to sum up with the original entries to $U$.

To illustrate, consider $3$-way arrays of format $L\times M\times N$
(already augmented with slack variables). Working
over alphabet $\{0,\dots,u\}$, define upper bounds $u_{i,j,k}:=u$
for original message variables and $u_{i,j,k}:=U$ for slack variable.
Then the set of possible messages that the receiver has to choose from is

\begin{equation}\label{feasible4}
\F\ :=\ \{x\in\Z^{L\times M\times N} \ :\ \sum_i x_{i,j,k}
= \sum_j x_{i,j,k}=\sum_k x_{i,j,k}=U \,,\ 0\leq x_{i,j,k}\leq u_{i,j,k}\,\}\ .
\end{equation}
Choosing $L$ and $M$ to be relatively small and fixed,
$\F$ is again the set of integer points in an $N$-fold system.
Corollary \ref{separable_and_distance} part 2 now enables
the efficient solution of the {\em error-correcting decoding} problem
$$\min\{\|{\hat x}-{\bar x}\|_p\ :\ {\hat x}\in \F\}\ .$$

\bc{error-correction}
Fix $L,M$. Then $3$-way $l_p$ error-correcting decoding can be done in polynomial time.
\ec

%

\subsection{Hierarchically-constrained multi-way arrays} 
\label{Subsection: Hierarchically-constrained multi-way arrays}

The transportation and routing problem, as well as the error-correction problem,
have very broad and useful generalizations, to arrays of any dimension and to any
{\em hierarchical} sum constraints. We proceed to define such systems of arrays.

Consider $d$-way arrays $x=(x_{i_1,\dots,i_d})$ of size $M_1\times\cdots\times M_d$.
For any $d$-tuple $(i_1,\dots,i_d)$ with $i_j\in\{1,\dots,M_j\}\cup\{+\}$,
the corresponding {\em margin} $x_{i_1,\dots,i_d}$ is the sum of entries of
$x$ over all coordinates $j$ with $i_j=+$. The {\em support} of $(i_1,\dots,i_d)$
and of $x_{i_1,\dots,i_d}$ is the set $\supp(i_1,\dots,i_d):=\{j:i_j\neq +\}$
of non-summed coordinates. For instance, if $x$ is a $4\times5\times3\times2$
array then it has $12$ margins with support $H=\{1,3\}$ such as
$x_{3,+,2,+}=\sum_{i_2=1}^5\sum_{i_4=1}^2 x_{3,i_2,2,i_4}$.
Given a family $\H$ of subsets of $\{1,\dots,d\}$ and margin values
$v_{i_1,\dots,i_d}$ for all tuples with support in $\H$, consider the
set of integer nonnegative and suitably upper-bounded arrays with these margins,
$$\F_\H \ :=\  \left\{\,x\in\Z^{M_1\times\cdots \times M_d}
\ : \ x_{i_1,\dots,i_d}\,=\,v_{i_1,\dots, i_d} \,,\ \ \supp(i_1,\dots,i_d)\in\H
\,,\ \ 0\leq x_{i_1,\dots,i_d}\leq u_{i_1,\dots,i_d}\right\}\ .$$
The congestion-avoiding transportation problem over $\F_\H$ is to find
$x\in\F_\H$ minimizing a given separable convex cost
$\sum_{i_1,\dots,i_d} f_{i_1,\dots,i_d}(x_{i_1,\dots,i_d})$.
The error-correcting decoding problem over $\F_\H$ is to estimate an
original message as ${\hat x}\in\F_\H$ minimizing a suitable $l_p$-distance
$\|{\hat x}-{\bar x}\|_p$ to a received message $\bar x$.

Again, for long arrays, that is, of format $M_1\times\cdots\times M_{d-1}\times N$
with $d$ and $M_1,\dots, M_{d-1}$ fixed and only the length (number of layers) $N$
variable, the set $\F_\H$ is the set of feasible points in an $N$-fold systems and,
as a consequence of Corollary \ref{separable_and_distance}, we can solve
both problems in polynomial time.

\bc{hierarchical}
Fix any $d,M_1,\dots,M_{d-1}$ and family $\H$ of subsets of $\{1,\dots,d\}$.
Then congestion-avoiding transportation and error-correcting decoding
over $\F_\H$ can be solved in polynomial time for any array length
$M_d:=N$ and any margin values $v_{i_1,\dots,i_d}$ for all tuples
$(i_1,\dots,i_d)$ with support in $\H$.
\ec

\section{Proofs of Theorems \ref{Theorem: Linear greedy
augmentation algorithm} and \ref{Theorem: Convex greedy augmentation
algorithm}}
\label{Section: Proofs to main theorems}

In this section we finally prove Theorems \ref{Theorem: Linear greedy
augmentation algorithm} and \ref{Theorem: Convex greedy augmentation
algorithm}. For this, we employ the following fact.

\begin{lemma}[Theorem 3.1 in Ahuja et al. \cite{Ahuja+Magnanti+Orlin}]
\label{Lemma: Sufficient improvement leads to polynomial algorithm}

Let $H$ be the difference between maximum and minimum objective
function values of an (integer valued) optimization problem.

Suppose that $f^k$ is the objective function value of some solution of
a minimization problem at the k-th interation of an algorithm and
$f^*$ is the minimum objective function value. Furthermore, suppose
that the algorithm guarantees that for every iteration $k$,
\[
(f^k-f^{k+1})\geq\beta(f^k-f^*)
\]
(i.e., the improvement at iteration $k+1$ is at least $\beta$ times
the total possible improvement) for some constant $0<\beta<1$ (which
is independent of the problem data). Then the algorithm terminates
in $O((\log H)/\beta)$.
\end{lemma}

\subsection{Proof of Theorem \ref{Theorem: Linear greedy
augmentation algorithm}}\label{Subsection: Proof of linear greedy algorithm}

Let $\Delta$ denote the least common multiple of all non-vanishing
maximal subdeterminants of $A$. Note that the encoding length
$\log\Delta$ is polynomially bounded in the encoding lengths of the
input data $A$, $u$, $b$ and $c$. Hence, the objective function values
of two vertices are either the same or differ by at least $1/\Delta$.

Let $f^0=\Delta\cdot c^\intercal z^0$ denote the normalized objective value
of the initially given feasible solution and by $f^1,f^2,\ldots$
denote the normalized objective values of the vertices
$z^1,z^2,\ldots$ that we reach at the end of the second steps of the
augmentation procedure. Note that the difference $H$ between maximum
and minimum normalized objective function values of $\LP_{A,u,b,f}$
has an encoding length $\log H$ that is polynomially bounded in the
encoding lengths of the input data $A$, $u$, $b$ and $c$. We now show
that
\[
(f^k-f^{k+1})\geq\beta(f^k-f^*)
\]
holds for $0<\beta=1/n<1$ and conclude by Lemma \ref{Lemma:
Sufficient improvement leads to polynomial algorithm}, that we only
have to enumerate $O((\log H)n)$, that is polynomially many, vertices.

Cosider the vector $z^*-z^k\in\ker(A)$. There is some orthant
$\Orthant_j$ such that $z^*-z^k\in\ker(A)\cap\Orthant_j$. Hence, we
can write
\[
z^*-z^k=\sum_{i=1}^n \alpha_i g_i
\]
for some $\alpha_i\in\R_+$ and $g_i\in\GLP(A)\cap\Orthant_j$,
$i=1,\ldots,n$. As $\alpha_i g_i$ has the same sign pattern as
$z^k-z^*$, one can easily check that the components of $z^k+\alpha_i
g_i$ lie between the components of $z^k$ and of $z^*$. Hence they are
nonnegative. As $g_i\in\ker(A)$, we have $Ag_i=0$ and thus
$A(z^k+\alpha_i g_i)=Az^k=b$ for any choice of
$i=1,\ldots,n$. Consequently, $z^k+\alpha_i g_i$ is a feasible
solution for any choice of $i=1,\ldots,n$. Finally, we have
\[
\Delta\cdot c^\intercal(z^k-z^*)=\sum_{i=1}^n \Delta\cdot c^\intercal(-\alpha_i g_i)=\sum_{i=1}^n \Delta\cdot c^\intercal(z^k-(z^k+\alpha_i g_i))
\]
from which we conclude that there is some index $i_0$ such that
\[
\Delta\cdot c^\intercal(z^k-(z^k+\alpha_{i_0} g_{i_0}))=\Delta\cdot c^\intercal(-\alpha_{i_0} g_{i_0})\geq\frac{1}{n}\sum_{i=1}^n \Delta\cdot c^\intercal(-\alpha_i g_i)=\frac{1}{n}\Delta\cdot c^\intercal(z^k-z^*)=\frac{1}{n}(f^k-f^*).
\]
Note that a greedy choice for an augmentating vector cannot make a
smaller augmentation step than the vector $\alpha_{i_0} g_{i_0}$. Thus,
\[
f^k-f^{k+1}\geq \Delta\cdot c^\intercal(z^k-(z^k+\alpha_{i_0} g_{i_0}))\geq\frac{1}{n}(f^k-f^*).
\]
This proves Part (a).

The proof of Part (b) is is nearly literally the same. Clearly, in the
integer situation, we may choose $\Delta=1$. If $z^1,z^2,\ldots$
denote the vectors that we reach from our initial feasible solution
$z^0$ via greedy augmentation steps, we only have to be careful about
the choice of $\beta$. In the integer situation, we
need to choose $\beta=1/(2n-2)$, since for the integer vector
$z^*-z^k\in\ker(A)\cap\Orthant_j$ at most $2n-2$ vectors from
the Hilbert basis of $C_j=\ker(A)\cap\Orthant_j$ are needed to
represent each lattice point in $C_j\cap\Z^n$ as a nonnegative
\emph{integer} linear combination of elements in
$\GIP(A)\cap\Orthant_j$ \cite{Seboe:90}. Thus, we need to apply
$O((\log H)(2n-2))=O((\log H)n)$ augmentation steps, a number being
polynomial in the encoding length. \eoproof

\subsection{Proof of Theorem \ref{Theorem: Convex greedy
augmentation algorithm}}\label{Subsection: Proof of convex greedy algorithm}


In \cite{Hemmecke:Z-convex,Murota+Saito+Weismantel}, it was shown that
$\GIP(A,C)$ allows a representation
\[
  (z^*-z^k,-C(z^*-z^k))=\sum_{i=1}^{2(n+s)-2} \alpha_i (g_i,-Cg_i),
\]
where each $\alpha_i\in\Z_+$ and where each $(g_i,-Cg_i)$ lies in the
same orthant as $(z^*-z^k,-C(z^*-z^k))$. It follows again from the
results in \cite{Seboe:90} that at most $2(n+s)-2$ summands are
needed. Similarly to the proof of Theorem \ref{Theorem: Linear greedy
augmentation algorithm}, we can already conclude from this
representation that $z^k+\alpha_i g_i$ is feasible for all
$i=1,\ldots,2(n+s)-2$.

Moreover, in \cite{Murota+Saito+Weismantel} it was shown that for such
a representation superadditivity holds, that is,
\[
\bar{f}(z^*)-\bar{f}(z^k)\geq\sum_{i=1}^{2(n+s)-2} [\bar{f}(z^k+\alpha_i g_i)-\bar{f}(z^k)]
\]
and thus, rewritten,
\[
\bar{f}^k-\bar{f}^*=\bar{f}(z^k)-f(z^*)\leq\sum_{i=1}^{2(n+s)-2} [\bar{f}(z^k)-\bar{f}(z^k+\alpha_i g_i)].
\]
Therefore, there is some index $i_0$ such that
\[
\bar{f}^{k}-\bar{f}^{k+1}=\bar{f}(z^k)-\bar{f}(z^k+\alpha_{i_0} g_{i_0})\geq
\frac{1}{2(n+s)-2} [\bar{f}(z^k)-\bar{f}(z^*)]=\frac{1}{2(n+s)-2}(\bar{f}^k-\bar{f}^*),
\]
and the result follows from Lemma \ref{Lemma: Sufficient improvement
leads to polynomial algorithm}. \eoproof

\end{document}